\documentclass[12pt]{article}
\usepackage{amsfonts,amsthm,amsmath,amssymb}
\newtheorem{Theorem}{Theorem}[section]
\newtheorem{Lemma}[Theorem]{Lemma}
\theoremstyle{definition}

\newtheorem{Remark}[Theorem]{Remark}
\newtheorem*{Proof}{Proof}
\newcommand{\boxtensor}{{\Box\kern-9.03pt\raise1.42pt\hbox{$\times$}}}

\numberwithin{equation}{section}
\newcounter{elno}                % This to number lists

\newcounter{example}[section]

\newcommand{\Spec}{{\mathrm{ Spec }}}

\title{ The Non-properness of the functor of $F$-trivial bundles}

\author{V.B. Mehta\footnote{This paper was presented by the first 
author at a conference for Peter Russell at McGill University, Montreal 
in June 2009. He would like to thank the organizers D. Daigle, R. 
Ganong, J. Hurtubise, M. Koras and S. Lu for the invitation and 
hospitality.} ~and~ S. Subramanian\footnote{Vikram Mehta 
passed away on 4$^{\mathrm{th}}$ June, 2014.}} 
\begin{document}
\maketitle

\begin{abstract}
We study the properness of the functor of $F$-trivial bundles by relating it 
to the base change question for the fundamental group scheme of Nori. 
\end{abstract}

\section{Introduction}

Let $X$ be a non-singular projective variety over an algebraically 
closed field  
of arbitrary characteristic with a very ample line bundle $H$ on $X$.
The notion of a coherent torsion-free sheaf being stable or semistable 
with respect to $H$ is now classical \cite{N3,L1}. In particular their 
moduli spaces with fixed Chern classes have been constructed [loc. 
cit]

In particular the property of the functor of semistable sheaves being 
proper is of crucial importance for  $ \dim X =1$, this was proved 
by 
Seshadri and then in general by Langton \cite{L2}.  For 
chi-semistability, or Giesekei- Maruyama semistability, this was proved 
by Mehta-Ramanathan \cite{MR} and Maruyama. The properness of the 
semi-stable  functor for $G$-bundles was also considered by Ramanathan 
for curves  in characteristic zero, (see also Balaji-Seshadri \cite{BS} and Faltings \cite{F}) 
then by Balaji-Parameswaran \cite{BP} for curves in characteristic $p$ and then by 
Heinloth \cite{H2} and Gomez-Langer-Sols-Schmitt \cite{GLSS} for arbitrary varieties in 
characteristic $p$.

Denninger-Werner consider a  slightly more general question \cite{DW}. Let 
$k$ be an algebraically closed  field  of characteristic $p$ and 
$W=W(k)$ 
the ring of Witt vectors over $k$, with function field $K$ and residue 
field $k$. Let $X\rightarrow  \ \nolinebreak \Spec(K)$ be a smooth, projective
absolutely irreducible curve over $K$.  Assume that $V$ is a semistable 
bundle of degree zero over $X$.  

They ask the following question:

Does there exist a model $\bar{X}$ of $X$, and an extension $\bar{V}$
of $V$ to $\bar{X}$, such that for each irreducible component $Y_i$ of 
$\bar{X}_k$, the restriction of $V$ to the normalization of 
$(Y_i)_{{\mathrm{red}}}$ 
%has the following properties  $V/Nor[(Y_i)_{red}]$ 
is strongly  semistable.  Note that the original $V$ on $X_K$ may be 
considered strongly semistable as characteristic $K=0$.  Such a bundle 
$V$ on $X_K$ is said to be have strong semistable reduction.  Suppose 
there exists a finite  morphism  $f:Z\rightarrow X$ over $K$ and a model 
$\bar{Z}$ of $Z$ such that for every irreducible component $Z_i$ of 
$\bar{Z}_k$,  the bundle $f^{*}(V)$ is strongly semistable on each 
normalization of $Z_{i_{\mathrm{red}}}$. Then $V$ is said to have potentially 
strong semistable reduction. Note the analogies with the semistable 
reduction of vector bundles and principal bundles mentioned earlier.

For vector bundles with a strong semistable reduction, Denninger-Werner 
show that there are functorial isomorphisms of ``parallel transport'' 
along etale paths between the fibres of $V_{\bar{K}} $ on $X_{\bar{K}}$, 
where ${\bar{K}}$ is the algebraic closure of $K$. See also Hackstein 
\cite{H1} for a similar discussion on $G$-bundles.

In another direction, Madhav Nori had introduced the fundamental 
group scheme of a reduced projective scheme $X$ over $k$, denoted by 
$\Pi^N(X)$, \cite{N1, N2}.  This is defined by assigning a Tannaka group 
to the Tannaka category of essentially finite vector bundles on 
$X$ [loc.cit].  In these papers, Nori had made 2 conjectures:

\medskip

(1) If $X$ and $Y$ reduced, complete schemes  over $k$, then 
$\Pi^N(X\times_k Y)\simeq \Pi^N(X)\times_k \Pi^N(Y)$.

\medskip

(2) If $l$ is an algebraically closed field extension of $k$, then the 
canonical map $\Pi^N(X_l) \rightarrow \Pi^N(X) \otimes_k l$ is an isomorphism.

\medskip

The present authors had proved conjecture (1) in \cite{MS1}, using the the 
notion of an ``$F$-trivial vector bundle''. They had also introduced the 
local fundamental group-scheme of $X$ denoted by $\Pi^{\mathrm{loc}}(X)$, using 
the 
Tannaka category of $F$-trivial bundles on $X$ \cite{MS2}.  In [loc. cit]  
they had also proved some necessary and sufficient conditions for the 
second conjecture of Nori to be valid.  In an attempt to prove the 
second conjecture of Nori, they had formulated the following question:

\paragraph*{Question 1:}  Let $X$ be a smooth projective curve, and let 
$S_0$ be a smooth affine curve with a smooth completion $S$.  Let $V_0$ 
be a vector bundle on $X \times S_0$ such that for every $s$ in $S_0, \ 
V\mid X \times \{s\}$ is $F$- trivial on $X$.  Then can $V_0$ be 
extended to a vector bundle $V$ on $X \times S$ such that for every  $s$ 
in $S, V\mid X \times \{s\}$ is $F$-trivial on $X$?

This may be thought of as a properness theorem for $F$-trivial vector 
bundles on $X$.

We may also consider the following question:

\paragraph*{Question 2:} Let $X$ be a nonsingular projective curve on 
$k$ and let $T_0$, be any smooth affine curve.  Let $V$ be a vector 
bundle on $X \times T_0$ such that
\begin{enumerate}
\item for all $t\in T_0$, the bundle 
$V\mid X\times \{t\}$ is $F$-trivial on $X$.
\item  for all $t$ in a non-empty open subset $U$ of $T_0$, the bundle  
$V\mid X \times \{t\}$ is stable on $X$ (and also $F$-trivial on $X 
\times \{t\}$) for all $t\in U$.
\end{enumerate}

Then is the classifying map $c: T_0\rightarrow U_X(r,0)$ constant?

Here $U_X(r,0)$ denotes the moduli spaces of rank $r$ and degree $0$ 
semistable vector bundles on $X$.  Note that any $F$-trivial vector 
bundle on $X$ is strongly semistable  of degree 0 \cite{MS2}.

In this paper we prove that an affirmative answer to Question 1 leads  
to an affirmative answer to Question 2.  It is important to note here 
that an affirmative answer to Question 2 would prove Nori's second 
conjecture.  In fact, Nori's second conjecture is equivalent to Question 
2. [Section 3].

But Christian Pauly has given a counter-example to Nori's second 
conjecture \cite{P}.  He constructs a nonconstant family of stable, 
$F$-trivial vector bundles, which is not constant. Therefore, Question 2, 
is false, hence Question 1 has also a negative answer. 

This also shows that in the equicharacteristic $p$ case, the question of 
Denninger-Werner also has a negative answer when one fixes a smooth 
and projective model for $X$, that is when the special fibre is a 
smooth, 
projective curve. For all the results used  here, about stability, 
semistability, $F$-trivial bundles and the precise  statements of Nori's 
2 conjectures, we refer to \cite{N3, MS2}.

\section{Formulation and some Lemmas}

Here we collect some basic facts about $F$-trivial bundles.  Throughout 
we work over an algebraically closed field $k$ of characteristic $p$.  If $X$ 
is a scheme (reduced) of finite type over $k$, we denote by $F$, 
the Frobenius map $X\rightarrow X$.  As $k$ is assumed to be 
perfect, we do not distinguish between the geometric Frobenius and the 
absolute Frobenius.  We have 

\paragraph*{Definition:} A vector bundle $V$ on $X$ is said to be 
$F$-trivial if $F^*(V) \simeq {\cal O}_X^{\oplus r}$, where $r$ =rank $V$ 
\cite{MS1}.

\begin{Remark}  For any integer $n >1$, we could define on 
$F^n$-trivial bundles on $X$ as a bundle $V$ on $X$ such that 
$F^{n*}(V) \simeq {\cal O}_X^{\oplus r}, r$=rank $V$.  But for ease of 
notation we assume $ n=1$.  See however, the remarks at the end of 
section 3.  We rephrase Questions 1 and 2 as Statements:
\end{Remark}

\paragraph*{Statement 1} Let $X$ be a nonsingular projective curve and 
$S_0$ a smooth affine curve, with smooth completion $S$(everything is 
defined over $k$).  Let $V_0$ be a vector bundle on $X \times S_0$ such 
that for all $s \in S_0, \ V_s: = V\mid X\times \{s \}$ is 
$F$-trivial  on $X$.  Then $V_0$ can be extended to a bundle $V$ on $X
\times S$ such that for all $s \in S, \ V_s := V\mid X \times \{ s \}$ 
is $F$-trivial on $X$.
%\end{Theorem}

Note that if such an extension exists, then it is unique, as $F$-trivial 
bundles are semistable, and by Langton's Theorem. Using Statement 1, we 
shall prove

\paragraph*{Statement 2}  Let $X$ be a nonsingular projective curve and 
$T_0$ a smooth affine curve.  Assume that there exists a vector bundle 
$V_0$ on $X \times T_0$ such that

\medskip

1) for all $t \in T_0, V_t := V_0 \mid X \times \{ t\}$ is $F$-trivial 
on $X$.

\medskip

2) for all $t$ in a non-empty open subset $U$ of $T_0$, $V_t$ is stable on 
$X$. Then the family $V_0$ is constant, that is the classifying map 
$c_{T_0}: T_0 \rightarrow U_X(r,0)$ is constant, where $U_X (r,0)$ is the 
moduli space of rank r and degree 0 semistable bundles on $X$. 
%\end{Theorem}

\bigskip

Assuming Statement 1, we prove Statement 2 in a sequence of Lemmas: 

\medskip

Let $V_0$ on $X \times T_0$ be as in Statement 2.  By Statement 1, $V_0$ can be 
extended to a vector bundle $V$ on $X \times T$, where $T$ is a smooth 
completion of $T_0$, such that for all $ t \in T, \ V_t : =V \mid X 
\times \{ t \}$ is $F$-trivial on $X$.  Then we have 

\begin{Lemma}
For any $x \in X$, consider the bundle $V_x :=  V\mid \{ x \} \times \ T$ on 
$T$. Then $\{ V _x\}, \ x \in X$, considered as a family of bundles on 
$T$ parameterized by $X$, is constant. That is, $V_x \simeq V_y$ as 
bundles on $T$, for any pair of points $x, y$ in $X$.
\end{Lemma}

\begin{Proof}
Let $F_X : X \rightarrow X$ be the Frobenius map of $X$.  Consider $(F_X 
\times \mathrm{id}_T) ^* (V)$ on $X \times T$. As $F_X^*(V_t)$ is trivial on 
$X$, for any $ t \in T$, we have $(F\times \mathrm{id}_T)^*(V) \simeq 
p_2^*(W)$ for some vector bundle on $W$ on $T$, by semicontinuity.  
Hence $(F_X \times \mathrm{id}_T)^*(V) \mid \{ x \} \times T \simeq W$ on $T$.  
But clearly, $V\mid \{ x \} \times T$ and $(F_X \times \mathrm{id}_T)^*(V) \mid 
\{ x \} \times T$ are isomorphic as bundles on $T$,  as 
$F_X: X\rightarrow X$ is surjective .  Hence the Lemma.  
\end{Proof}

In what follows we shall call $W$ the parameter bundle on $T$.

\begin{Lemma}
Without loss of generality, $W$ may be assumed to have degree $0$ on $T$.
\end{Lemma}

\begin{Proof}
If degree $W=0$, we are through. Otherwise, let degree $W=d$ and rank $W 
=r$ = rank $V$.  We can certainly find a map. $f:Z \rightarrow T$, where 
$Z$ is a smooth curve and a line bundle $L$ on $Z$ such that $r\deg L 
+ \deg f^*(W) = 0$  Consider the family $(\mathrm{id}_X\times f)^*(V) \otimes p^*L^{-1}$ on
$X\times Z$.  The parameter bundle for this family is clearly
$f^*(W)\otimes L^{-1}$, which has degree 0 on $Z$.  The new family on
$X \times Z$ has the same properties as $V$ on $X \times T$.  And if
$c_Z:Z\rightarrow U_X(r,0)$ is constant, so is  the map $c_T:
T\rightarrow U_X(r,0)$ where $c_T$ and $c_Z$ are the classifying map for
$V_0$ on $X\times T$ and $(\mathrm{id}_X \times f)^*(V) \otimes  p^*_2 L^{-1}$ on
$X\times Z$ respectively.  So we may assume that the degree of the
parameter bundle $W = 0$.
\end{Proof}

Now we make an assumption that will be removed in Section 3 
\begin{enumerate}
\item[{$(*):$}] 
The ground field $k$ is the algebraic closure of $F_p$, that is 
$k \simeq \bar{F}_p$. 
\end{enumerate}
With this assumption we have

\begin{Lemma}
$W$ is strongly semistable on $T$.
\end{Lemma}

\begin{Proof} If not, $F^{m^*}(W)$ has a strong  Harder-Narasimhan 
filtration: $ 0 =A_0 \subset A_1 \subset \cdots \subset A_n =W$ with each 
$A_i/A_{i-1}$ strongly semistable and $\mu (A_i/A_{i-1}) > \mu (A_j/A_{j-1})$ if 
$ i < j$.  Denote each $A_i/A_{i-1}$ by $B_i$.  Just as in Lemma~2.3, 
$\exists$ $f: Z_i\rightarrow T$ such that $f^*(B_i)$ has degree 0, with 
$Z_i$ a smooth curve.  But $f^*(B_i)$ is still strongly semistable on 
$Z_i$.  
\iffalse
Now consider $F^{a^*} (f^* B_i),\  a=1,2,\ldots$.  This is a 
family of semistable bundles on $Z_i$, of fixed rank $r$ and degree 0.  
As we are working over $k =\bar{F}_p$, there are only finitely many rational
points of the moduli space of vector bundles of degree $0$ over a finite 
field. Hence given a strongly semistable bundle $B$ of degree $0$, defined 
over a finite field, the set $\{F^{m^*} B,\ m=1,2,\ldots\}$ is a finite set. Hence 
there exists $a,b, \ a\neq b$ with 
$F^{a^*}(f^* B_i)\simeq F^{b^*} (f^* B_i)$.  
\fi
This implies that 
$f^*(B_i)$ is essentially finite on $Z$ \cite{S}. This further implies 
that there exists a smooth projective curve $S$ and a map $g: S 
\rightarrow Z$ with $g^*(f^* B_i)$ trivial on $S$.  Trivializing  
$B_1,\ldots, B_n$ this way, we get finally a smooth projective curve $R$ 
and a map $h: R \rightarrow T$ such that $h^*(W)$ is a direct sum,
$$
h^*(W) \simeq L^{\oplus r_1}_1 \oplus L^{\oplus r_2}_2 \cdots \oplus L^{\oplus r_s}_s
$$
for some line bundles $L_i, 1\leq i \leq s$ on $R$ with $\deg L_1> \deg 
L_2 \ldots >\deg L_s$ and some positive 
integers $r_i$, with $\sum r_i =r$. Consider the family $(\mathrm{id}_X \times 
h)^*(V)$ on $X \times R$.  If the classifying map $c_R : R \rightarrow 
U_X(r,0)$ is constant, then the classifying map  $c_T: T \rightarrow U_X 
(r,0)$ is also constant.  Hence we may assume, that on $T$ itself the 
parameter bundle  is a  direct sum, 
$ W\simeq L^{\oplus r_1}_1\oplus \cdots \oplus L_s^{\oplus r_s}$ with each $L_i \in \ \mathrm{Pic}(T)$.

Now consider $V$ as a family of vector bundles on $T$ parametrized by 
$X$.  Consider the relative Harder-Narasimhan  filtration of $V$ on $T$, 
relative to 
$X$.  By semi-continuity $\exists$ vector bundles $V_1,\ldots, V_s$ on 
$X$ of ranks $ r_1, \ldots, r_s$, and a filtration on $X \times T$
$$
0=W_0 \subset W_1 \cdots \subset W_s =V \eqno{(1)}
$$
such that each $W_i/W_{i-1} \simeq p_1^* (V_i) \otimes p_2^*(L_i)$
Now consider the family $(F_X \times \mathrm{id}_T)^* (V)$ on $X \times T$.  This 
bundle is isomorphic to $p^*_2(W)$ by Lemma~2.2. So 
$$
p_2^*(W) \simeq p_2^* L_1^{\oplus r_1}\oplus \cdots \oplus p^*_2 L^{\oplus r_s}_s 
\eqno{(2)}
$$
Apply $(F_x \times \mathrm{id}_T)^* $ to filtration (1), we get
$$
0\subset F^*_X(W_1)\ldots \subset F^*_X (W_s) =F^*_X(V) \eqno{(3)}
$$
with each $F^*_X(W_i)/ F^*_X(W_{i-1}) \simeq p^*_1 (F^*_X V_i) \otimes 
p^*_2(L_i)$.

Compare filtration (2) and (3):

It is clear that $p^*_2 L^{\oplus r_1}_1$ has  no maps to $p^*_1 (F_X^* 
V_i) \otimes p^*_2 (L_i)$ for $i>1$. as degree $L_1 > $ degree $L_2 
\ldots >$ degree $L_s$.

Hence $p^*_2 L^{\oplus r_1}_1$ injects into $p_1^*(F^*_X V_1)\otimes 
p_2^*(L_1)$ on $X \times T$.  Tensoring by $p^*_2 (L^{-1}_1)$, we get 
an injection of the trivial bundle of rank $r$, into $p^*_1(F^*_X 
V_1)$ on $X \times T$.  Hence on $X$, we get an injection of the trivial 
bundle into $F^*_X V_1$.  This implies that degree $F^*_X V_1 \geq 
0$ hence degree $V_1 \geq 0$. Choose any closed point $ t \in T$ and 
restrict filtration (1) to $X \times \{ t \}$ we get an injection of 
$V_1$ inside $V_t$.  But degree $V_1= 0$ and $V_t$ is stable for a a 
general $t \in T$, which is a contradiction.  Hence $W$ is strongly 
semistable on $T$.

\end{Proof}

\begin{Theorem}
The classifying map $c_T: T\rightarrow U_X (r,0)$ is constant.
\end{Theorem}

\begin{Proof}
Consider the sequence of bundles on $T$, given by $F^{n^*} (W), \ 
n=1,2,\ldots $.  They are all semistable of degree $0$. As we are 
working over $k=\bar{F}_p$, we must have $F^{n^*} (W) \simeq F^{m^*} 
(W)$ for some positive integers $m,n$ with $m\neq n$.  By \cite{S}, $W$ is 
essentially finite on $T$, that is there exists a smooth projective curve 
$Z$ and a  map $h: Z \rightarrow T$ such that $h^*(W)$ is trivial.  
Consider the family $(\mathrm{id}_X\times h)^*(V)$ on $X \times Z$.  This has 
parameter bundle $h^*(W)$, which is trivial.  So $(\mathrm{id}_X \times h)^*(V) 
\simeq p_1^*(A)$ for some vector bundle $A$ on $X$.  But then the 
classifying map $c_Z: Z \rightarrow U_X(r,0)$ is constant, so the classifying map $c_T: T \rightarrow U_X(r,0)$ is also constant.
\end{Proof}

\section{Main result}

Now we remove the assumption $(*)$ in Section 2.  So let $k$ be any 
algebraically closed field of characteristics $p$ and  $V_0$ on $X 
\times T_0$ be as in Statement~2.  By Statement~1,  $V_0$ extends to a 
family $V$ of $F$-trivial bundles on $X$, parameterized by $T$,  where $T$ is a 
smooth completion of $T_0$  Then we have

\begin{Theorem} The classifying map $c_T: T\rightarrow 
U_X(r,0)$ is constant.
\end{Theorem}

\begin{Proof}  
$W$ is defined by the isomorphism $(F_X \times \mathrm{id}_T)^*(V) \simeq p_2^*(W)$ on $X \times T$.  We may assume that there exists an algebra $R$, 
finitely generated over $F_p$ such that for $V$ on $X \times T$,  there 
exist models:
$$
\begin{array}{llll} 
1) & X_R \rightarrow R & \mbox{for} &   X\\[1mm]
2) & T_R \rightarrow R & \mbox{for}  & T.\\[1mm]
3) & V_R \rightarrow X_R  &  \mbox{for}& V\\
4) & W_R \rightarrow T_R & \mbox{for} & W.
\end{array}
$$
We may also assume that there exist open 
nonempty subsets $O_1$ and $O_2$ of $\Spec (R)$ such that 1) and 2) below are
satisfied:

\medskip

1) For  every geometric point $\Spec (\omega) \rightarrow O_1$ with image a 
closed point $m$ of $O_1$ , the bundle $V_\omega$ is a family of $F$-trivial 
bundles on $X_\omega$, parameterized by $T_\omega$.  This is seen as 
follows : the bundle $W$ on $T$ has a model $W_R \rightarrow T_R$ and the 
isomorphism $(F_X \times \mathrm{id}_T)^* (V) \simeq  p^*_2(W)$ can be spread out 
over $O_1$.  This proves that for every closed point $t$ of $T_\omega$, 
the bundle $V_t$ is $F$-trivial on $X_\omega\times \{t\}$.

\medskip

2)We may also assume that there exists a nonempty open subset $O_2$ of 
$\Spec R$ such that for every geometric point $\Spec (\omega) \rightarrow O_2$ 
with image a closed point $m$ of $O_2$, the family $V_\omega$ on 
$X_\omega \times T_\omega$ is generically stable.  This is proved as  
follows:

Let $U_{X_R} (r.0) \rightarrow \ \Spec(R)$  be the relative moduli space 
of semistable bundles of rank $r$ and degree 0 on the fibers of 
$X_R\rightarrow R$ \cite{L1}. Let $U^s_{X_R}$ be the open submoduli space 
of stable bundles.  If $K$ is the quotient field of $R$, then we have 
assumed that the image of $c_K:T_K \rightarrow U_{X_K}$ intersects 
$U_{X_K}^s$.  So there exists a nonempty open subset $O_2$ of $\Spec (R)$ 
such that $c_m:T_m \rightarrow U_{X_m}$  intersects $U^s_{X_m}$ for all 
closed points $m$ in $O_2$.  So we may assume, without loss of 
generality, that $O_1 =O_2 =\Spec R$.  For any closed point $m$ in 
$\Spec R$ we know that $c_m : T_m \rightarrow U_{X_m}$ is a constant map.  It 
easily follows that $c_K: T_K \rightarrow U_{X_K}$ is a constant map, 
thus finishing the proof that Statement  1 implies Statement 2 over arbitrary 
algebraically closed fields in characteristic~$p$.
\end{Proof}

\begin{Remark} For any integer $n >1$, we may define an $F^n$-trivial bundle $V$ on $X:$ we ask that $F^{n^*} (V)={\mathcal O}^{\oplus 
r}_X , r=$ rank $V$.  It is trivial to check that  the proof that 
Statement 1 implies Statement 2 goes through without any changes.
\end{Remark}

\begin{Remark}  But Statement 2 is false!! More precisely, 
Christian Pauly \cite{P} has produced a non-constant  family of stable 
bundles, trivialized by $F^4$, the fourth power of Frobenius on a curve 
of genus 2 in characteristic 2 . He has 
done this by a careful study of the  Verschiebung map $V: U_X(r,0) 
\rightarrow U_X(r,0)$, induced by the Frobenius map $F: X \rightarrow 
X$.  This also shows that the answer to the question of Denninger-Werner 
is also negative, if one works over a fixed smooth model of $X$.  In 
fact the second conjecture of Nori implies Statement~1. We provide a sketch 
proof in the following:
\end{Remark}

\begin{Theorem}  The second conjecture  of Nori implies the 
properness of the functor of $F$-trivial bundles (i.e., Statement 1).
\end{Theorem}

\begin{Proof}  Let $V_0$ over $X \times T_0 $ be as in Statement 1, and 
let $T$ be a smooth completion of $T_0$.  We may assume, without loss of 
generality, that $T-T_0$ is one point, say $\infty$.  We first extend $V_0$ 
to a family of semistable  bundles on $X\times T$,  
that is, $V_\infty$ is semistable.  Note that $V_\infty$ is not unique, 
but $\mathrm{gr}(V_\infty)$ is unique.  If $V_t$ is stable for some point $t\in 
T_0$, then $V_t$ is stable for all $t$ in $U$, where $U$ is open in 
$T_0$. By \cite{MS2}, we know that the set of isomorphism classes 
$V_t$ is finite.  This implies that the map $c_U : U\rightarrow 
U_X(r,0)$ is constant and $c_U(t)$ is a fixed bundle $V$ in $U^s_X(r,0)$ for all $t$ in $U$.  As $T$ is connected, we must have $c_U(t)=V$ in $U^s_X(r,0)$ for 
all $t$ in $T$, in particular $V_\infty$ is $F$-trivial on $X$.  
Assume that $V_t$ is strictly semistable for all $t \in T_0$.  Let 
$[A_{t,i}]$ be the set of stable components of $V_t, \ t $ in $T_0$, and 
$1\leq i < r =$ rank $V$.  Again by \cite{MS2}, the set of isomorphism 
classes of $\{A_{t,i}\}$ is finite.  But this implies that the map 
$c_{T_0}\rightarrow U_X(r,0)$ has a finite image, hence constant as 
$T_0$ is connected.  So again, we get
$\mathrm{gr}(V_\infty)= \mathrm{gr}(V_t)$  for all $t$ in $T$.

But this implies that $V_\infty$ is strongly semistable. (All the 
Frobenius pull back are semistable of degree 0).  Now consider the 
family $F^*_X(V_t), t \in T$.  This a family of trivial bundles 
converging to a semistable  bundle $F^*(V_\infty)$.  So 
$F_X^*(V_\infty)$ is also trivial, so $V_\infty$ is $F$-trivial.
\end{Proof}

\begin{Remark}  To complete the circle of ideas, we note that 
Statement 2 implies the second conjecture of Nori.  We prove this now. 
\end{Remark}

\begin{Theorem}  The validity of Statement 2 implies the second 
conjecture of Nori.  
\end{Theorem}

\begin{Proof}  Let $k \subset K$ be algebraically closed fields of 
characteristic $p$ and let $V$ be a stable  $F$-trivial bundle on 
$X_K:= X {\otimes_k} K$.  We check the criterion in \cite{MS2}.  So we have to show that 
$V$ is  defined over $k$, that is over $X$.  We can find $R$, a finitely 
generated algebra with quotient field $L \subset K$ such that $X_K$ has a 
model $X_R \rightarrow R$ and $V$ has a model $V_R \rightarrow X_R$.  By 
cutting down $\Spec (R)$ suitably, we may assume that for all geometric 
points $\Spec (\omega)\rightarrow \Spec(R)$,  the bundles $V_\omega$ are 
stable and $F$ trivial on $X_\omega$. For any curve $T_0$ in $\Spec (R)$ the 
family $V_R\mid X_R  \times_R T_0$ is constant by assumption. So the map \linebreak
$c_R: \Spec(R) \rightarrow U_X(r,0)$ is constant, where  $r$ =rank $V$.  
But this means that $V$ is defined over $k$, that is, $V$ comes from $X$.
\end{Proof}

\begin{Remark}
It maybe true that a generalized version of the properness theorem is
true for the Nori fundamental group scheme. This would also imply the  
generalized conjecture of Denninger-Werner. We hope to return to these 
questions in the future.
\end{Remark}

\noindent S.Subramanian \\
School of Mathematics \\
TIFR,Navy Nagar \\
Mumbai 400005 \\ 
email:subramnn@math.tifr.res.in \\

\end{document}